\documentclass[12pt]{amsart}
\usepackage{graphicx,amssymb,latexsym}

\theoremstyle{plain}
\newtheorem{theorem}{Theorem}[section]

\newtheorem{lemma}[theorem]{Lemma}
\newtheorem{proposition}[theorem]{Proposition}
\newtheorem{definition}[theorem]{Definition}

\theoremstyle{remark}
\newtheorem{remark}{Remark}[section]
\newtheorem{example}[remark]{Example}

\newcommand{\ol}{\overline{\lambda}}
\newcommand{\om}{\overline{\mu}}
\newcommand{\ul}{\underline{\lambda}}
\newcommand{\um}{\underline{\mu}}
\newcommand{\lsm}{\lambda/\mu}
\begin{document}

\title[Schur and skew Schur functions]
{Equality of Schur and skew Schur functions}%
\author{Stephanie van Willigenburg}%
\address{Department of Mathematics, University
of British Columbia, Vancouver, BC  V6T 1Z2}%
\email{steph@math.ubc.ca}%

\thanks{The author was supported in part by the
National Sciences and Engineering Research Council of Canada.}%

\begin{abstract}
We determine the precise conditions under which any skew Schur
function is equal to a Schur function over both infinitely and
finitely many variables.
\end{abstract}
\maketitle
\section{Introduction}\label{background}

Littlewood-Richardson coefficients arise in a variety of contexts.
The first of these is that they are the structure constants in the
algebra of symmetric functions with respect to the basis of Schur
functions. Another instance is as the multiplicities of
irreducible representations in the tensor product of
representations of the symmetric group. A third occurrence is as
intersection numbers in the Schubert Calculus on a Grassmanian.
Thus knowing their values has an impact on a number of fields. In
this paper we calculate when certain coefficients are $0$ or $1$
by determining when a skew Schur function is equal to a Schur
function. Although multiplicity free products have been studied in
~\cite{Stembridge} our determination will reveal more precisely
when certain coefficients are $0$ and when they are $1$. The
related question of when two ribbon Schur functions are equal has
been answered recently in \cite{Billera_&_co}, which revealed many
new equalities of Littlewood-Richardson coefficients.

The remainder of this note is structured as follows. In the rest
of Section ~\ref{background} we review the definitions required.
This is followed by the two main theorems in which we give
straightforward conditions that prescribe when a skew Schur
function is equal to a Schur function over both infinitely and
finitely many variables.

\subsection{Schur and skew Schur functions}

We say that a list of positive integers $\lambda =\lambda
_1\geq\lambda _2 \geq\ldots \geq\lambda _k$ whose sum is $n$ is a
\emph{partition} of $n$, denoted $\lambda \vdash n$. We call the
$\lambda _i$ the \emph{parts} of $\lambda$. A partition with at
most one part size is called a \emph{rectangle} and a partition
with exactly two different part sizes is called a \emph{fat hook}.
If $\lambda = \lambda _1\lambda _2 \ldots \lambda _k\vdash n$ then
we define the \emph{(Ferrers) diagram} $D_\lambda$ to be the array
of left justified boxes with $\lambda _i$ boxes in the $i$-th row
for $1\leq i\leq k$. If we transpose $D_\lambda$ we obtain another
diagram $D_{\lambda '}$ known as the \emph{conjugate} of
$D_\lambda$ and we refer to $\lambda$ and $\lambda '$ as
\emph{conjugate partitions}. Furthermore, for any column $c$ in
$D_\lambda$ we denote by $l(c)$ the number of boxes in $c$ and
refer to $l(c)$ as the \emph{length} of $c$. Where the context is
clear we abuse notation and refer to $D_\lambda$ as $\lambda$.

\begin{example}
\begin{center}
$D_{4322}=$\begin{picture}(100,50)(-30,-5)
\put(0,0){\line(1,0){20}} \put(0,10){\line(1,0){20}}
\put(0,20){\line(1,0){30}} \put(0,30){\line(1,0){40}}
\put(0,40){\line(1,0){40}}

\put(0,40){\line(0,-1){40}} \put(10,40){\line(0,-1){40}}
\put(20,40){\line(0,-1){40}} \put(30,40){\line(0,-1){20}}
\put(40,40){\line(0,-1){10}}
\end{picture}
\end{center}
\end{example}

We define a \emph{(Young) tableau} $T$ of \emph{shape} $\lambda$
to be a filling of the boxes of $D _\lambda$ with positive
integers. If the filling is such that the integers in each row
weakly increase, whilst the integers in each column strictly
increase we say that $T$ is a \emph{semi-standard} tableau.
\begin{example}\hfill$\begin{array}{cccc} 1&1&2&3\\
2&3&4&\\
4&4&&\\
5&7&& \end{array}$ is a semi-standard tableau of shape $4322$.
\end{example}

If $\mu=\mu _1\ldots \mu _k\vdash m$, $\lambda = \lambda _1\ldots
\lambda _l\vdash n$ where $k\leq l$ and $m\leq n$ such that $\mu
_i\leq \lambda _i$ for $1\leq i\leq k$ then we define the
\emph{skew diagram} $D_{\lambda /\mu}$ to be the array of boxes
that appear in $D_\lambda$ but not in $D_\mu$. For our purposes
$D_{\lambda/\mu}$ will always be \emph{connected}, that is to say,
for any pair of adjacent rows in $D_{\lambda/\mu}$ there exists at
least one column in which they both have a box.

\begin{example}
\begin{center}
$D_{4322} =$ \begin{picture}(100,50)(-30,-5)
\put(0,0){\line(1,0){20}} \put(0,10){\line(1,0){20}}
\put(0,20){\line(1,0){30}} \put(0,30){\line(1,0){40}}
\put(0,40){\line(1,0){40}}

\put(0,40){\line(0,-1){40}}
\put(10,40){\line(0,-1){40}}
\put(20,40){\line(0,-1){40}}
\put(30,40){\line(0,-1){20}}
\put(40,40){\line(0,-1){10}}

\end{picture}\hfill
$D_{21}= $\begin{picture}(100,50)(-30,-5)
\put(0,0){\line(1,0){10}} \put(0,10){\line(1,0){20}}
\put(0,20){\line(1,0){20}}

\put(0,20){\line(0,-1){20}}
\put(10,20){\line(0,-1){20}}
\put(20,20){\line(0,-1){10}}
\end{picture}
\hfill $D_{4322/21}=$
\begin{picture}(100,50)(-30,-5)
\put(0,0){\line(1,0){20}}
\put(0,10){\line(1,0){20}}
\put(0,20){\line(1,0){30}}
\put(10,30){\line(1,0){30}}
\put(20,40){\line(1,0){20}}

\put(0,20){\line(0,-1){20}} \put(10,30){\line(0,-1){30}}
\put(20,40){\line(0,-1){40}} \put(30,40){\line(0,-1){20}}
\put(40,40){\line(0,-1){10}}
\end{picture}
\end{center}
\end{example}

Again when the context is clear we refer to $D_{\lambda /\mu}$ as
$\lambda /\mu$. Similarly we define \emph{skew tableaux} and
\emph{semi-standard skew tableaux} by appropriately inserting the
adjective skew in the above definitions for tableaux and
semi-standard tableaux.

\begin{definition}
Given a (skew) tableau $T$ and a set of variables $x_1, x_2, \ldots$ we
define the monomial $x^T$ to be
$$x^T:= x_1^{t_1}x_2^{t_2}\ldots$$
where $t_i$ is the number of times $i$ appears in $T$. Let
$\lambda, \mu$ be partitions such that $\lambda/\mu$ is a (skew)
diagram, then we define the corresponding \emph{(skew) Schur
function} $s_{\lambda/\mu}$ to be
$$s_{\lambda/\mu}=\sum x^T$$
where the sum is over all semi-standard (skew) tableaux $T$ of
shape $\lambda/\mu$.
\end{definition}

The set of all Schur functions $s_\lambda$ (i.e. $s_{\lambda
/\mu}$ where $\mu=\emptyset$) forms a basis for the algebra of
symmetric functions, $\Lambda$, which is a subalgebra of
$\mathbb{Q}[x_1, x_2, \ldots ]$.

Since it can be easily shown that skew Schur functions are
symmetric it follows that skew Schur functions can be written as a
linear combination of Schur functions. To be more precise we need
to recall two more notions: that of the reading word and the
content of a (skew) tableau. Firstly, given a (skew) tableau, $T$,
we say its \emph{reading word}, $w(T)$, is the entries of the
tableau read from top to bottom and \emph{right} to \emph{left}.
Given a reading word we say it is \emph{lattice} if as we read it
from left to right the number of $i$'s we have read is at least as
large as the number of $i+1$'s we have read e.g. $1213$ is
lattice, however, $1132$ is not as when we have read $113$ the
number of $3$'s we have read is greater than the number of $2$'s.
Secondly, the \emph{content} of a (skew) tableau, $c(T)$, is a
list $t_1t_2t_3\ldots$ where, as before, $t_i$ is the number of
times $i$ appears in $T$.

We are now ready to express any skew Schur function as a linear
combination of Schur functions.

\begin{proposition}~\cite[A1.3.3]{ECII}
Let $\lambda, \mu, \nu$ be partitions such that $\lambda/\mu$ is a
(skew) diagram then
$$s_{\lambda/\mu}=\sum _\nu c^\lambda _{\mu\nu} s_\nu$$
where $c^\lambda _{\mu\nu}$ is the number of semi-standard (skew)
tableaux $T$ such that
\begin{enumerate}
\item the shape of $T$ is $\lambda /\mu$
\item $c(T)=\nu$
\item $w(T)$ is lattice.
\end{enumerate}\end{proposition}

\begin{example}
$$s_{32/1}=s_{31}+s_{22}.$$
\end{example}

\begin{remark}
The $c ^\lambda _{\mu\nu}$ are known as Littlewood-Richardson
coefficients, and the above method of computing them is known as
the Littlewood-Richardson rule. There are many other methods for
computing the $c^\lambda _{\mu\nu}$ such as Zelevinsky's pictures
 or Remmel and Whitney's reverse numbering, and the interested reader may
 wish to consult, say, ~\cite{Fulton} for further details. However, it is the
Littlewood-Richardson rule that will allow us to determine our
results most succinctly.
\end{remark}

\section{Equality of Schur and skew Schur functions}

Before we state our main result let us define an involution on
diagrams. Given a diagram $\lambda$ let $\lambda ^\circ$ be the
(skew) diagram that is the diagram $\lambda$ rotated by $180
^\circ$.

\begin{example}
\begin{center}
$4322 =$ \begin{picture}(100,50)(-30,-5) \put(0,0){\line(1,0){20}}
\put(0,10){\line(1,0){20}} \put(0,20){\line(1,0){30}}
\put(0,30){\line(1,0){40}} \put(0,40){\line(1,0){40}}

\put(0,40){\line(0,-1){40}} \put(10,40){\line(0,-1){40}}
\put(20,40){\line(0,-1){40}} \put(30,40){\line(0,-1){20}}
\put(40,40){\line(0,-1){10}}
\end{picture}
$4322 ^\circ =$  \begin{picture}(100,50)(-30,-5)
\put(0,0){\line(1,0){40}} \put(0,10){\line(1,0){40}}
\put(10,20){\line(1,0){30}} \put(20,30){\line(1,0){20}}
\put(20,40){\line(1,0){20}}

\put(0,10){\line(0,-1){10}}
\put(10,20){\line(0,-1){20}}
\put(20,40){\line(0,-1){40}}
\put(30,40){\line(0,-1){40}}
\put(40,40){\line(0,-1){40}}
\end{picture}\end{center}
\end{example}

\begin{theorem}\label{skew_equal_schur}
For  partitions $\lambda, \mu,\nu$
$$s_{\lambda/\mu}=s_\nu \mbox{ if and only if } \lambda/\mu = \nu \mbox{ or } \nu ^\circ .$$
\end{theorem}

\begin{proof}
The reverse implication follows by Exercise 7.56(a) ~\cite{ECII},
which yields that
$$s_\nu=s_{\nu^\circ}.$$

For the forward implication we need only show that if $\lambda
/\mu$ is not $\nu$ or $\nu ^\circ$ for some diagram $\nu$ then
$\lambda/\mu$ has more than one filling whose reading word is
lattice.

Consider the skew tableau $T$ of shape $\lambda/\mu$ where each column
$c$ is filled with the integers $1,\ldots , l(c)$ in increasing
order. This filling is clearly lattice. Now since $\lambda/\mu$ is
not a diagram $\nu$ nor a (skew) diagram $\nu ^\circ$ where $\nu$
is a diagram, consider the first row $i$ where $\lambda /\mu$
fails to be either $\nu$ or $\nu ^\circ$ for some diagram $\nu$
(i.e. if $\lsm$ is truncated at row $i-1$ then we obtain a
(rotated) diagram, but this is no longer true if $\lsm$ is
truncated at row $i$). Moving from right to left note the first
entry $j$ in $T$ which does not have $i-1$ entries above it. Form
the reading word of $T$ upto this entry and note the smallest
integer $i\geq k>j$ for which the number of occurrences of $k$ is
strictly less than the number of occurrences of $k-1$. Change $j$
to $k$ and change all entries below it in that column by adding
$k-j$ to the existing entry to form a new skew tableau $T'$ of shape
$\lambda/\mu$. Since $w(T')$ is clearly lattice, we are done.
\end{proof}

\section{Equality and $GL(n)$ or $SL(n)$ characters}

The set of all Schur functions restricted to the variables $x_1,
\ldots , x_n$, obtained by setting $x_m=0$ for $m>n$, forms a
basis for the algebra of symmetric polynomials $\Lambda _n$. Skew
Schur functions in $\Lambda _n$ can be expressed in terms of the
Schur functions by
$$s_{\lambda /\mu}(x_1, \ldots ,x_n)= \sum _\nu c^\lambda_{\mu\nu} s_\nu(x_1, \ldots ,x_n)$$
and thus we can ask when $s_{\lambda /\mu}(x_1, \ldots
,x_n)=s_\nu(x_1, \ldots ,x_n).$ In terms of representation theory
this yields when certain multiplicities in the tensor products of
irreducible representations of $SL(n, \mathbb{C})$ (or polynomial
representations of $GL(n, \mathbb{C})$) will be 0 and when they
will be 1.

Clearly if $s_{\lambda/\mu}= s_\nu$ in $\Lambda$ then the result
holds in $\Lambda _n$, however the converse may not be true as an
$s_{\lambda/\mu}$ comprising of a sum of $s_\nu$ only one of which
has less than $n+1$ parts could exist. However, the search for
such an $s_{\lambda/\mu}$ is greatly reduced as the converse may
not be true only when the length of the longest column in
$\lambda/\mu$ is equal to $n$ by

\begin{lemma}\label{n_variables}
Let $\lambda, \mu$ be partitions such that the length of the
longest column in $\lambda /\mu$ is $m$.
\begin{enumerate}
\item If $m>n$ then $s_{\lambda/\mu}=0$ in $\Lambda _n$.
\item If $m<n$ then $s_{\lambda/\mu}= s_\nu$ in $\Lambda _n$ if
and only if $s_{\lambda/\mu}= s_\nu$ in $\Lambda$.
\end{enumerate}\end{lemma}

\begin{proof}
The first result is immediate from the definitions. The reverse
direction of the second result has already been discussed, thus it
only remains to show that if $m<n$ and $s_{\lambda/\mu}\neq s_\nu$
in $\Lambda$ then $s_{\lambda/\mu}\neq s_\nu$ in $\Lambda _n$.

Consider the skew tableau $T$ of shape $\lambda/\mu$ where each column
$c$ is filled with the integers $1,\ldots , l(c)$ in increasing
order. Since $s_{\lambda /\mu} \neq s_\nu$ in $\Lambda$ this
implies $\lambda /\mu$ is not a diagram $\nu$ nor a (skew) diagram
$\nu ^\circ$, so as in the proof of Theorem
~\ref{skew_equal_schur} consider the first row $i$ where $\lambda
/\mu$ fails to be $\nu$ or $\nu ^\circ$ for some diagram $\nu$.
Moving from right to left note the first entry $j$ in column $c$
of $T$ which does not have $i-1$ entries above it, form the
reading word upto this entry, $w^\prime$, and note the smallest
integer $i\geq k>j$ for which the number of occurrences of $k$ is
strictly less than the number of occurrences of $k-1$. If
$l(c)+k-j\leq n$ then change the $j$ to $k$ and change all entries
below it in $c$ by adding $k-j$ to the existing entry to form a
new skew tableau $T'$ of shape $\lambda /\mu$. If not then find the
largest entry in $w^\prime$, $k^\prime <n$,  and fill the
$n-k^\prime$ lowest boxes in $c$ with $k^\prime +1, \ldots ,n$ to
form $T'$. Since $w(T')$ is clearly lattice, the result follows.
\end{proof}

\begin{example}
If $n=5$ then the following semi-standard skew tableaux illustrate 
$T'$ in the situation $l(c)+k-j\leq n$ and $l(c)+k-j> n$ respectively.
\begin{center}
$\begin{array}{cccc}
&&1&1\\
&2&2&\\
1&3&&\\
2&4&&
\end{array} \qquad
\begin{array}{ccc}
&&1\\
&&2\\
&&3\\
&1&4\\
&2&\\
1&5&\end{array}$
\end{center}
\end{example}

Thus, from here on we shall assume that the length of the longest
column in $\lambda /\mu$ is $n$. Before we reveal the analogous
result to Theorem ~\ref{skew_equal_schur} let us define two
operations on (skew) diagrams.

\begin{definition}
Let $\lambda = \lambda _1 \ldots \lambda _k$ and $\mu =\mu_1\ldots
\mu _l$ be partitions such that $\lambda /\mu$ is a (skew)
diagram. Let $c$ be a column of longest length in $\lambda /\mu$
and $\ol$ and $\om$ be partitions such that
\begin{enumerate}
\item $\ol ' = (\lambda _1 '+r)(\lambda _2 '+r)(\lambda _3
'+r)\ldots (\lambda _c '+r)\lambda _{c+1} '\ldots \lambda
_{\lambda _1} '$

$\om ' = (\mu _1 '+r)(\mu _2 '+r)(\mu _3 '+r)\ldots (\mu _c
'+r)\mu _{c+1} '\ldots \mu _{\mu _1} '$ or
\item $\ol ' = (\lambda _1 '+r)(\lambda _2 '+r)(\lambda _3
'+r)\ldots (\lambda _{c-1} '+r)\lambda _{c} '\ldots \lambda
_{\lambda _1} '$

$\om ' = (\mu _1 '+r)(\mu _2 '+r)(\mu _3 '+r)\ldots (\mu _{c-1}
'+r)\mu _{c} '\ldots \mu _{\mu _1} '$
\end{enumerate}
and $\lambda/\mu$ is a (skew) diagram then we say $\ol/\om$ is a
\emph{shearing} of $\lambda/\mu$.
\end{definition}

\begin{remark}
Intuitively we can interpret this definition as creating a diagram
$\ol/\om$ from $\lambda/\mu$ by choosing a column of longest
length and sliding it and every column to the left of it down $r$
boxes, or sliding it and every column to the right of it up $r$
boxes.
\end{remark}

\begin{example}
The first two skew diagrams are shearings of $442$ whilst the
third is not.

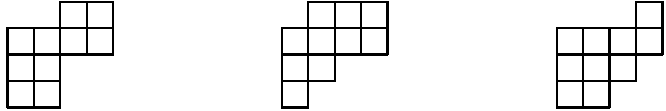
\begin{picture}(100,50)(-30,-5)
\put(0,0){\line(1,0){20}} \put(0,10){\line(1,0){20}}
\put(0,20){\line(1,0){40}} \put(0,30){\line(1,0){40}}
\put(20,40){\line(1,0){20}}

\put(0,30){\line(0,-1){30}} \put(10,30){\line(0,-1){30}}
\put(20,40){\line(0,-1){40}} \put(30,40){\line(0,-1){20}}
\put(40,40){\line(0,-1){20}}
\end{picture}
\begin{picture}(100,50)(-30,-5)
\put(0,0){\line(1,0){10}} \put(0,10){\line(1,0){20}}
\put(0,20){\line(1,0){40}} \put(0,30){\line(1,0){40}}
\put(10,40){\line(1,0){30}}

\put(0,30){\line(0,-1){30}} \put(10,40){\line(0,-1){40}}
\put(20,40){\line(0,-1){30}} \put(30,40){\line(0,-1){20}}
\put(40,40){\line(0,-1){20}}
\end{picture}
\begin{picture}(100,50)(-30,-5)
\put(0,0){\line(1,0){20}} \put(0,10){\line(1,0){30}}
\put(0,20){\line(1,0){40}} \put(0,30){\line(1,0){40}}
\put(30,40){\line(1,0){10}}

\put(0,30){\line(0,-1){30}} \put(10,30){\line(0,-1){30}}
\put(20,30){\line(0,-1){30}} \put(30,40){\line(0,-1){30}}
\put(40,40){\line(0,-1){20}}
\end{picture}
\end{example}

\begin{definition}\label{fattening}
Let $\ul$ and $\um$ be partitions such that
\begin{enumerate}
\item $\ul ' = (a+b)^ic^{j+k}$ and $\um '= b^{i+k}$ if $a\neq c$
or
\item $\ul ' = (a+b)^i(\nu _1+b)\ldots (\nu _k+b)c^{j}$ and $\um '=
b^{i+k}$, where $\nu=\nu _1\ldots \nu _k$ is a partition, if $a=c$
\end{enumerate}
and $\ul/\um$ is a (skew) diagram then we say $\ul/\um$ is a
\emph{fattening} of $(a^ic^j)'$ and $(\ul/\um)^\circ$ is a
\emph{fattening} of $((a^ic^j)')^\circ$.
\end{definition}

\begin{remark}
Intuitively we can interpret this definition as creating a diagram
say $\ul/\um$ from a fat hook or rectangle in the following way.
If we have a fat hook $(a^ic^j)'$ then we shear the rightmost
column of length $a$ and all the columns to the left of it down by
$b$ boxes. We then insert a rectangle $k^{c-b}$ such that the
result is a (skew) diagram. If we have a rectangle $(a^{i+j})'$
then we shear a column and all the columns to the left of it down
by $b$ boxes. We then insert a diagram such that the result is a
(skew) diagram. A similar interpretation follows for
$(\ul/\um)^\circ$.
\end{remark}

\begin{example}
All three skew diagrams are fattenings of $444$.

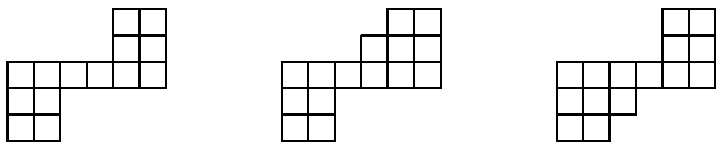
\begin{picture}(100,70)(-30,-5)
\put(0,0){\line(1,0){20}} \put(0,10){\line(1,0){20}}
\put(0,20){\line(1,0){60}} \put(0,30){\line(1,0){60}}
\put(40,40){\line(1,0){20}}\put(40,50){\line(1,0){20}}

\put(0,30){\line(0,-1){30}} \put(10,30){\line(0,-1){30}}
\put(20,30){\line(0,-1){30}} \put(30,30){\line(0,-1){10}}
\put(40,50){\line(0,-1){30}}\put(50,50){\line(0,-1){30}}
\put(60,50){\line(0,-1){30}}
\end{picture}
\begin{picture}(100,70)(-30,-5)
\put(0,0){\line(1,0){20}} \put(0,10){\line(1,0){20}}
\put(0,20){\line(1,0){60}} \put(0,30){\line(1,0){60}}
\put(30,40){\line(1,0){30}}\put(40,50){\line(1,0){20}}

\put(0,30){\line(0,-1){30}} \put(10,30){\line(0,-1){30}}
\put(20,30){\line(0,-1){30}} \put(30,40){\line(0,-1){20}}
\put(40,50){\line(0,-1){30}}\put(50,50){\line(0,-1){30}}
\put(60,50){\line(0,-1){30}}
\end{picture}
\begin{picture}(100,70)(-30,-5)
\put(0,0){\line(1,0){20}} \put(0,10){\line(1,0){30}}
\put(0,20){\line(1,0){60}} \put(0,30){\line(1,0){60}}
\put(40,40){\line(1,0){20}}\put(40,50){\line(1,0){20}}

\put(0,30){\line(0,-1){30}} \put(10,30){\line(0,-1){30}}
\put(20,30){\line(0,-1){30}} \put(30,30){\line(0,-1){20}}
\put(40,50){\line(0,-1){30}}\put(50,50){\line(0,-1){30}}
\put(60,50){\line(0,-1){30}}
\end{picture}
\end{example}

For convenience we extend the notion of fattening to all (skew)
diagrams by defining the fattening of $\lambda /\mu$ to be
$\lambda /\mu$ if $\lambda /\mu$ is any (skew) diagram other than
those referred to in Definition ~\ref{fattening}. In addition, for
clarity of exposition, we denote by $\widetilde{\lambda/\mu}$ any
(skew) diagram that has been derived from $\lambda/\mu$ via some
combination of shearings or fattenings.

\begin{theorem}
For partitions $\lambda, \mu, \nu, \eta =\eta _1\eta_2 \ldots$
where $\eta _i$ is the number of columns of length at least $i$ in
$\lsm$
$$s_{\lambda/\mu}=s_\eta \mbox{ in $\Lambda _n$ if and only if } \lambda/\mu =
\widetilde{\nu} \mbox{ or } \widetilde{\nu ^\circ}.$$
\end{theorem}

\begin{proof}
For the reverse implication it is straightforward to check that if
$\lambda/\mu = \widetilde{\nu} \mbox{ or } \widetilde{\nu ^\circ}$
then the only semi-standard skew tableau $T$ of shape $\lambda/\mu$ whose reading
word is lattice has each column $c$ filled with the integers $1,
\ldots ,l(c)$ in increasing order.

The forward implication will follow once we show that if
$\lambda/\mu \neq \widetilde{\nu} \mbox{ or } \widetilde{\nu
^\circ}$ for some diagram $\nu$ then $\lambda/\mu$ has more than
one filling utilising the integers $1, \ldots ,n$ whose reading
word is lattice.

Consider the skew tableau $T$ of such a shape $\lambda/\mu$ whose
reading word is lattice where each column $c$ is filled with the
integers $1, \ldots ,l(c)$ in increasing order. If $\lsm$ has a
column $c_1$ such that $l(c_1)<n$ and $c_1$ contains at least one
box with no box to the right of it, and a column $c_2$ to the
right of $c_1$ such that $l(c_2)<n$ then if $l(c_1)\leq l(c_2)$
change the entry in the last box of $c_1$ to $l(c_2)+1$ otherwise,
unless the box at the head of $c_1$ and the column immediately to
the right of it are in the same row, change the entry $l(c_2)$ in
$c_1$ to $l(c_2)+1$ and increase all entries below it in $c_1$ by
1 to form a new skew tableau $T'$ of shape $\lsm$ whose reading word is
lattice.

Thus if $\lsm$ does not satisfy these criteria then $\lsm$ must be
of the form
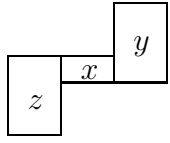
\begin{center}
\begin{picture}(100,70)(-30,-5) \put(0,0){\line(1,0){20}}
\put(20,20){\line(1,0){40}} \put(0,30){\line(1,0){40}}
\put(40,50){\line(1,0){20}}

\put(0,30){\line(0,-1){30}}  \put(20,30){\line(0,-1){30}}
\put(40,50){\line(0,-1){30}} \put(60,50){\line(0,-1){30}}

\put(7.5, 10){$z$}\put(27.5,21.25){$x$}\put(47.5,32.5){$y$}
\end{picture}\end{center}
where $x$ is a rectangle, and by Lemma ~\ref{n_variables} and what
we have already proved $y$ (and $z$) must consist of a
non-rectangular (skew) diagram $\delta$ or $\delta ^\circ$, for
some diagram $\delta$, whose column lengths are all less than $n$
with a non-negative number of columns to the right or left of it
that consist of $n$ boxes. If $c_s$ is the column of shortest
length $l(c_s)$ in $y$ then change the entry $l(c_s)$ in the first
column of $x$ to $l(c_s)+1$ and increase all the entries below it
in that column by 1. If the first column of $x$ does not contain
the entry $l(c_s)$ then change the entry in the last box to
$l(c_s)+1$. In both instances we form a new skew tableau $T'$ of shape
$\lsm$ whose reading word is lattice.

Hence the number of columns in $x$ must be zero and $\lsm$ must be
of the form
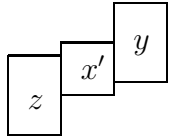
\begin{center}
\begin{picture}(100,70)(-30,-5)
\put(0,0){\line(1,0){20}} \put(0,30){\line(1,0){20}}
\put(20,15){\line(1,0){20}}\put(20,35){\line(1,0){20}}
\put(40,20){\line(1,0){20}}\put(40,50){\line(1,0){20}}

\put(0,30){\line(0,-1){30}}  \put(20,35){\line(0,-1){35}}
\put(40,50){\line(0,-1){35}} \put(60,50){\line(0,-1){30}}

\put(7.5, 10){$z$}\put(27.5,21.25){$x'$}\put(47.5,32.5){$y$}
\end{picture}\end{center}
where the length of every column of $x'$ is $n$ and $y$ (and $z$)
consists of a non-rectangular (skew) diagram $\delta$ or $\delta
^\circ$, for some diagram $\delta$, whose column lengths are all
less than $n$ with a non-negative number of columns  to the right
(respectively left) that consist of $n$ boxes. For clarity of
exposition we identify $y$ and $z$ with the sub (skew) diagram
$\delta$ or $\delta ^\circ$ that they contain. Let $\delta$ and
$\varepsilon$ be diagrams and $c_l$, $c_s$ be the columns of
longest or shortest length in $\delta$ respectively.

If $y=\delta$ and $z=\varepsilon$ then it follows that
$x'$ must contain zero columns otherwise we can change the entries in the
first column of $\varepsilon$; every column in $\varepsilon$ must
be at least as long as the longest column of $\delta$; and the box
at the head of the leftmost column of $\delta$ and the rightmost
column of $\varepsilon$ are not in the same row. Change the entry
$l(c_l)$ in the first column of $\varepsilon$  to $l(c_l)+1$ and
increase all the entries below it in that column by 1.

If $y=\delta ^\circ$ and $z=\varepsilon ^\circ$ then $x'$ must
contain zero columns; every column in $\varepsilon ^\circ$ must be
no longer than the shortest column in $\delta ^\circ$; and the box
at the base of the leftmost column of $\delta ^\circ$ and the
rightmost column of $\varepsilon ^\circ$ are not in the same row.
Change the last entry in the first column of $\varepsilon ^\circ$
to $l(c_l)+1$.

If $y=\delta$ and $z=\varepsilon ^\circ$ then either $x'$ must
contain zero columns or the box at the base of the rightmost
column of $\varepsilon ^\circ$ is in the same row as the adjacent
column of length $n$. In either case we can apply an argument
similar to that above to either change the entries of the
rightmost column of $\varepsilon ^\circ$ from $l(c_s)$ downwards
by increasing them by 1 or change the entry in the last box to
$l(c_s)+1$.

Finally if $y=\delta ^\circ$ and $z=\varepsilon$ then the same
conditions must be satisfied as for the case $y=\delta$ and
$z=\varepsilon$. However, we can change the entries of the
rightmost column of $\varepsilon$ from $l(c_l)$ downwards by
increasing them by 1.

In each situation we have been able to create a new skew tableau $T'$
of shape $\lsm$ whose reading word is lattice, and having
eliminated all possibilities the result follows.
\end{proof}

\section*{Acknowledgements}

The author would like to thank Eric Babson for suggesting the
problem, Benjamin Young for his skew diagram drawing package and the
referee for their diligence. 


\end{document}